# CENTRAL LIMIT THEOREM FOR THE SOLUTION OF THE KAC EQUATION

BY ESTER GABETTA[1] AND EUGENIO REGAZZINI[2]

*Università degli Studi di Pavia and Università degli Studi di Pavia*

We prove that the solution of the Kac analogue of Boltzmann's equation can be viewed as a probability distribution of a sum of a random number of random variables. This fact allows us to study convergence to equilibrium by means of a few classical statements pertaining to the central limit theorem. In particular, a new proof of the convergence to the Maxwellian distribution is provided, with a rate information both under the sole hypothesis that the initial energy is finite and under the additional condition that the initial distribution has finite moment of order $2+\delta$ for some $\delta$ in $(0,1]$. Moreover, it is proved that finiteness of initial energy is necessary in order that the solution of Kac's equation can converge weakly. While this statement may seem to be intuitively clear, to our knowledge there is no proof of it as yet.

## 1. Introduction and presentation of new results.

1.1. *Introduction.* Marc Kac studied Boltzmann's derivation of a basic equation of kinetic theory by simplifying the problem to an $n$-particle system in one-dimension and, under suitable conditions, he got the following analogue of the Boltzmann equation:

$$
(1) \quad \begin{cases} \dfrac{\partial}{\partial t} f(v,t) = \dfrac{1}{2\pi} \displaystyle\int_{\mathbb{R}\times[0,2\pi)} \{f(v\cos\theta - w\sin\theta, t) \\ \qquad\qquad\qquad\qquad\qquad\quad \times f(v\sin\theta + w\cos\theta, t) \\ \qquad\qquad\qquad\qquad\qquad\qquad - f(v,t)f(w,t)\} \, dw\, d\theta, \\ f(v,0) = f_0(v) \qquad (t>0, v\in\mathbb{R}), \end{cases}
$$

Received October 2006; revised October 2007.

[1]Supported in part by Ministero dell'Istruzione, dell'Università e della Ricerca (MIUR Grant 2006/015821).

[2]Supported in part by Ministero dell'Istruzione, dell'Università e della Ricerca (MIUR Grant 2006/134526).

*AMS 2000 subject classifications.* 60F05, 82C40.

*Key words and phrases.* Central limit theorem, Kac's equation, Kologorov distance, Wild's sum.

where $f_0$ and $f(\cdot,t)$ denote the probability density functions of the velocity of each particle at time 0 and at time $t$, respectively. This problem admits a unique solution within the class of all probability density functions on $\mathbb{R}$. See, for example, Kac (1956), Kac (1959), McKean (1966), Cercignani (1975) and Diaconis and Saloff-Coste (2000).

Bobylëv (1984) proved that the Fourier transform $\phi(\cdot,t)$ of the solution $f(\cdot,t)$ of (1) must satisfy

$$
(2) \quad \begin{cases} \dfrac{\partial}{\partial t}\phi(\xi,t) = \dfrac{1}{2\pi}\int_0^{2\pi} \phi(\xi\sin\theta,t)\phi(\xi\cos\theta,t)\,d\theta - \phi(\xi,t), \\ \phi(\xi,0) = \phi_0(\xi) \quad (t>0, \xi\in\mathbb{R}), \end{cases}
$$

$\phi_0$ being the Fourier transform of $f_0$. Clearly, problem (2) is well defined for arbitrary (not necessarily absolutely continuous) probability measures $\mu(\cdot,t)$ and $\mu_0$ on the class $\mathscr{B}(\mathbb{R})$ of all Borel subsets of $\mathbb{R}$, provided that $\phi(\cdot,t)$ and $\phi_0$ are thought of as Fourier–Stieltjes transforms of $\mu(\cdot,t)$ and $\mu_0$, respectively.

The solution of (2)—which exists and is unique within the Fourier–Stieltjes transforms of all probability measures on $\mathscr{B}(\mathbb{R})$—can be expressed by means of the transform of the Wild series [see Wild (1951)], that is,

$$
(3) \quad \phi(\xi,t) = \sum_{n\geq 1} e^{-t}(1-e^{-t})^{n-1}\hat{q}_n^+(\xi;\phi_0) \qquad (t\geq 0, \xi\geq 0),
$$

where $\hat{q}_n^+$ can be found by recursion as

$$
\hat{q}_n^+(\xi;\phi_0) = \frac{1}{n-1}\sum_{j=1}^{n-1}\hat{q}_{n-j}^+(\xi;\phi_0)\circ\hat{q}_j^+(\xi;\phi_0) \qquad (n=2,3,\ldots),
$$

with $\hat{q}_1 := \phi_0$. The symbol $g_1 \circ g_2$, where $g_1$ and $g_2$ are characteristic functions, designates the *Wild product*

$$
g_1 \circ g_2(\xi) = \frac{1}{2\pi}\int_0^{2\pi} g_1(\xi\cos\theta)g_2(\xi\sin\theta)\,d\theta \qquad (\xi\in\mathbb{R}).
$$

Getting down to the approach to equilibrium of the solution of (1) as $t$ goes to infinity, according to Boltzmann's classical research, the entropy of $f(\cdot,t)$ should increase to its upper bound, $\log(\sigma\sqrt{2\pi e})$ with $\sigma^2 = \int_{\mathbb{R}} v^2 f_0(v)\,dv$, while $f$ tends to the Maxwellian function (viz., the Gaussian density with zero mean and variance $\sigma^2$)

$$
g_\sigma(v) = \frac{1}{\sigma\sqrt{2\pi}}e^{-v^2/(2\sigma^2)} \qquad (v\in\mathbb{R}).
$$

McKean (1966) argues that the Wild representation suggests a simpler explanation: the *central limit theorem for Maxwellian molecules*. With the aim of demonstrating the solidity of his argument, he starts by proving a new expression for $\hat{q}_n^+$, that is,

$$\hat{q}_n^+(\xi; \phi_0) = \sum_{\gamma \in G(n)} p_n(\gamma) c_\gamma(\xi; \phi_0), \tag{4}$$

where $c_\gamma$ denotes the *n-fold Wild product* of $\phi_0$ with itself performed according to an algebraic structure schematized by the element $\gamma$ of a class $G(n)$ of random trees with $n$ leaves. $p_n(\cdot)$ is a probability on the subsets of $G(n)$. See McKean (1967), Carlen, Carvalho and Gabetta (2000), Carlen, Carvalho and Gabetta (2005), Bassetti, Gabetta and Regazzini (2007). Then, considering the form of $c_\gamma$, with the aid of the Lindeberg version of the central limit theorem, McKean proves the following statement on the weak convergence of the probability distribution function $C_\gamma$, corresponding to $c_\gamma$, toward the Gaussian distribution function $G_\sigma(v) = \int_{-\infty}^v g_\sigma(x)\,dx$:

Set $\sigma^2 := \int_\mathbb{R} v^2 f_0(v)\,dv$ and let $\int_\mathbb{R} |v|^3 f_0(v)\,dv$ be finite. Then, for any $\delta > 0$, there are constants $c = c(\delta, f_0)$, $c_1 = c_1(\delta, f_0)$ and a positive integer $n_0 = n_0(\delta, f_0)$ such that

$$p_n\left(\left\{\gamma \in G(n) : \sup_{v \in \mathbb{R}} |C_\gamma(v) - G_\sigma(v)| > \delta\right\}\right) \leq c(\delta, f_0) n^{8/(3\pi)-1} \tag{5}$$

$$(n \geq n_0),$$

which leads to

$$\sup_{v \in \mathbb{R}} |F(v, t) - G_\sigma(v)| \to 0 \qquad (\text{as } t \to +\infty), \tag{6}$$

where $F(\cdot, t)$ denotes the probability distribution function which corresponds to the solution $\phi(\cdot, t)$ of (2).

1.2. *Presentation of new results.* The study of necessary and sufficient conditions under which (6) holds true, together with some hints to rate of convergence, is the main scope of the present paper. We will prove the following:

THEOREM 1. *Let $\mu_0$ be a nondegenerate probability measure on $\mathscr{B}(\mathbb{R})$ and let $F(\cdot, t)$ be the probability distribution function corresponding to the solution $\phi(\cdot, t)$ of the Kac equation (2). Then*

$$\sup_{v \in \mathbb{R}} |F(v, t) - G_\sigma(v)| \to 0 \qquad (\text{as } t \to +\infty)$$

*holds true if and only if $\sigma^2 := \int_\mathbb{R} x^2 \mu_0(dx)$ is finite.*

It is wellknown that (6) is valid when the initial energy is finite. See, for example, Carlen and Lu (2003). As far as the necessity of this condition is concerned, while it cannot be doubted from a physical intuitive standpoint, it seems that no proof of it has been advanced as yet. Moreover, our approach leads to state a rather precise quantitative evaluation of the rate of convergence. This result is contained in the next theorem, where $F_0$ and $F_{0,d}$ are probability distribution functions defined by

$$F_0(x) := \mu_0((-\infty, x]),$$
$$F_{0,d}(x) := \mu_0([-x, +\infty)) \qquad (x \in \mathbb{R}).$$

THEOREM 2. *Let $\mu_0$ be a nondegenerate probability measure on $\mathscr{B}(\mathbb{R})$, $\sigma^2 := \int x^2 \mu_0(dx)$ be finite and let $a$, $p$ be fixed numbers in $(0,1)$ and $(2,+\infty)$, respectively. Then, there is a strictly positive constant $A$ such that*

$$\sup_{x \in \mathbb{R}} |F(x,t) - G_\sigma(x)| \leq AM(t)^{1/5} + \tfrac{1}{2} \sup_{x \in \mathbb{R}} |F_0(x) - F_{0,d}(x)| e^{-t},$$

*where*

$$M(t) = \int_{|u| > \sigma(x_t)^{a-1}} u^2 \mu_0(du) \vee e^{-B_1 t} \vee e^{-B_2 t}$$

*for every $t \geq t_0 := \inf\{t : \int_{|u| > \sigma(x_t)^{a-1}} u^2 \mu_0(du) \leq 1\}$ with*

$$x_t := \exp\{-tcp\},$$

$$B_1 := acp, \qquad B_2 := 1 - 2\alpha_p - cp, \qquad c \in \left(0, \frac{1 - 2\alpha_p}{p}\right),$$

$$\alpha_p = \frac{1}{2\pi} \int_0^{2\pi} |\cos \theta|^p \, d\theta.$$

*Moreover, if $\overline{m}_{2+\delta} = \int |x|^{2+\delta} \mu_0(dx) < +\infty$ for some $\delta$ in $(0,1]$, then*

$$\sup_{x \in \mathbb{R}} |F(x,t) - G_\sigma(x)| \leq C_\delta \frac{\overline{m}_{2+\delta}}{\sigma^{2+\delta}} e^{-t(1-2\alpha_{2+\delta})} + \frac{1}{2} \sup_{x \in \mathbb{R}} |F_0(x) - F_{0,d}(x)| e^{-t},$$

*where $C_\delta$ is a universal constant (Berry–Esseen constant).*

Constant $A$ can be easily obtained by looking at the proof of Theorem 2 in Section 3.

The proofs of Theorems 1 and 2 rest on an idea which goes back to McKean (1966). In the present paper we go deep into that idea by providing a complete proof of the next basic theorem, in which $q_t(n) := e^{-t}(1-e^{-t})^{n-1}$, $n = 1, 2, \ldots$; $u^\infty$ is the probability measure on $\mathscr{B}([0, 2\pi)^\infty)$ which makes the coordinates of $[0, 2\pi)^\infty$ independent and uniformly distributed, and $\mu_0^\infty$ meets the same conditions with $[0, 2\pi)$ and $u$ replaced by $\mathbb{R}$ and $\mu_0$, respectively.

THEOREM 3. *For each $t > 0$, there are a probability space $(\Omega, \mathscr{F}, P_t)$ and random variables*

$$\tilde{\nu}_t : \Omega \to \mathbb{N},$$
$$\tilde{\gamma} : \Omega \to G := \bigcup_n G(n),$$
$$\tilde{\theta} := (\tilde{\theta}_1, \tilde{\theta}_2, \ldots) : \Omega \to [0, 2\pi)^\infty,$$
$$\tilde{x} := (\tilde{x}_1, \tilde{x}_2, \ldots) : \Omega \to \mathbb{R}^\infty$$

*with joint distribution*

$$P_t\{\tilde{\nu}_t = n, \tilde{\gamma} = \gamma, \tilde{\theta} \in A, \tilde{x} \in B\} = q_t(n)p_n(\gamma)\mathbb{1}_{G(n)}(\gamma)u^\infty(A)\mu_0^\infty(B)$$
$$(n \in \mathbb{N}, \gamma \in G, A \in \mathscr{B}([0, 2\pi)^\infty), B \in \mathscr{B}(\mathbb{R}^\infty))$$

*such that*

$$V_t := \sum_{j=1}^{\tilde{\nu}_t} \pi_j(\tilde{\gamma}, \tilde{\theta})\tilde{x}_j$$

*has probability distribution $\mu(\cdot, t)$, that is, the distribution corresponding to the solution $\phi(\cdot, t)$ of (2).*

Apart from the definition of functions $\pi_j$, that we postpone to Section 2, where a physical interpretation is given, Theorem 3 allows us to understand the connection between convergence to equilibrium of $\mu(\cdot, t)$ and central limit theorem: Indeed, $\mu(\cdot, t)$ is the distribution of $V_t$, that is, a sum of random variables. With respect to ordinary applications of the central limit theorem, here we have a *random number* ($\tilde{\nu}_t$) *of summands*, and these summands have (not stochastically independent) *random coefficients* ($\pi_j, \quad j = 1, \ldots$). But these difficulties can be avoided through a careful utilization of the features of the joint distribution of $(\tilde{\nu}_t, \tilde{\gamma}, \tilde{\theta}, \tilde{x})$. This way, one can provide complete proofs of Theorems 1 and 2 through simple adaptations of powerful classical results from probability theory. Actually, we are pursuing the object of tracing to the above very same kind of ideas the study of the trend to equilibrium (with rate information) under the most important (weak or strong) modes of convergence, both for the Kac model and for other models such as an "inelastic" version of (1)–(2) introduced in Pulvirenti and Toscani (2004), and the Boltzmann equation for Maxwellian molecules in case of spatially homogeneous initial data with uniform collision kernel [see, e.g., Carlen and Lu (2003)].

It is well to pause here and consider what will be involved in the arguments used to prove Theorems 1 and 2. First, we will provide a proof for Theorem 2 under the sole hypothesis that the initial energy is finite. It is apparent that this covers also the sufficiency part in Theorem 1. The line of reasoning, to obtain the rate of convergence in Theorem 2, consists in adapting the argument generally used in the proof of the classical Lindeberg–Feller version of the central limit theorem. As far as the necessity part in Theorem 1 is concerned—that is, convergence in distribution of $V_t$ implies that $\sigma^2$ is finite—we will resort to a method used in Fortini, Ladelli and Regazzini (1996) to prove central limit theorems for arrays of partially exchangeable random variables. The method rests on the fact that Theorem 3 entails conditional independence of the summands in the definition of $V_t$, given $(\tilde{\nu}_t, \tilde{\gamma}, \tilde{\theta})$. After denoting conditional distribution of $V_t$, given $(\tilde{\nu}_t, \tilde{\gamma}, \tilde{\theta})$, by $\Lambda_{\tilde{\nu}_t}$, the next step consists in proving that convergence in distribution of $V_t$, as $t \to +\infty$, implies that any increasing and diverging to infinity sequence of positive terms $t_1, t_2, \ldots$ contains a subsequence $(t_{n'})$ for which

(7)　the distribution of $\Lambda_{\tilde{\nu}_{t_{n'}}}$ weakly converges to the distribution of $\Lambda$,

$\Lambda$ being some (random) probability measure. Then, one combines (7) with the Skorokhod–Dudley representation to transform (7) into a statement about (almost sure) weak convergence of a suitably defined random distribution $\Lambda^*_{\tilde{\nu}^*_{t_{n'}}}$ toward a random probability measure $\Lambda^*$, where $\Lambda^*_{\tilde{\nu}^*_{t_{n'}}}$ has the distribution of $\Lambda_{\tilde{\nu}_{t_{n'}}}$, and $\Lambda^*$ has the distribution of $\Lambda$. At this stage, the central limit theorem is employed to deduce *necessary* conditions for the convergence of $\Lambda^*_{\tilde{\nu}^*_{t_{n'}}}$. Finally, one concludes by showing that these conditions boil down to the existence of a bounded variance for the initial distribution $\mu_0$.

As to organization of the rest of the paper, Section 2 includes, in addition to some necessary preliminary concepts and notation, a proof of Theorem 3. In Section 3 the reader can find the proofs of Theorems 1 and 2. The Appendix contains the proofs of a few preparatory propositions.

**2. Preliminaries and proof of Theorem 3.** The first part of the section contains elements necessary to the definition of the functions $\pi_j$ mentioned in Theorem 3. Recall that, if $\gamma$ is any McKean tree with $n \geq 2$ leaves, each node has either zero or two "children," a "left child" and a "right child" such as in Figure 1, where a few elements of $G(8)$ are visualized.

In each tree of $G(n)$ fix an order on the set of the $(n-1)$ nodes and, accordingly, associate the random variable $\tilde{\theta}_k$ with the $k$th node. See (a) in Figure 1. Moreover, call $1, 2, \ldots, n$ the $n$ leaves following a left to right order. See (b) in Figure 1. The number of generations which separate leaf $j$ from the "root" node is said to be the *depth* of $j$ (in symbols, $\tilde{\delta}_j$). With $\tilde{\delta}_{(1)}(\gamma)$ one denotes the *depth of the tree* $\gamma$, that is, $\min\{\tilde{\delta}_1(\gamma), \ldots, \tilde{\delta}_n(\gamma)\}$ if $\gamma \in G(n)$. The cardinality of $G(n)$ is the Catalan number $C_n = \binom{2n-2}{n-1}/n$; see Section 15 of Comtet (1970).

Now, for any leaf $j$ of $\gamma$ in $G(n)$, look at the path which connects $j$ and the "root" node at the top in ascending order. It consists of $\tilde{\delta}_j$ steps: the first one from $j$ to its "parent" node, the second from the "parent" to the "grandparent" of $j$, and so on. Define the product
$$\pi_j = \pi_j(\tilde{\gamma}, \tilde{\theta}) = \tau_1^{(j)} \cdots \tau_{\delta_j}^{(j)},$$
where $\tau_{\delta_j}^{(j)}$ equals $\cos \tilde{\theta}_k$ if $j$ is a "left child" or $\sin \tilde{\theta}_k$ if $j$ is a "right child" and $\tilde{\theta}_k$ is the element of $\tilde{\theta}$ associated to parent node of $j$; $\tau_{\delta_j - 1}^{(j)}$ equals $\cos \tilde{\theta}_m$ or $\sin \tilde{\theta}_m$ depending on if the "parent" of $j$ is, in its turn, a "left child" or a "right child," $\tilde{\theta}_m$ being the element of $\tilde{\theta}$ associated with the grandparent of $j$; and so on. For instance, as to leaf 1 in (a) of Figure 1, $\pi_1 = \cos \tilde{\theta}_4 \cdot \cos \tilde{\theta}_2 \cdot \cos \tilde{\theta}_1$ and, for leaf 6, $\pi_6 = \sin \tilde{\theta}_5 \cdot \cos \tilde{\theta}_3 \cdot \sin \tilde{\theta}_1$.

From the definition of $\pi_j$ one obtains
$$\sum_{j \in \gamma} \pi_j^2 = 1 \tag{8}$$
for every $\gamma$ in $G(n)$, with $n = 2, 3, \ldots$. It is worth extending (8) to $G(1)$ by setting $\pi_1 \equiv 1$ for the sole leaf of $\gamma$ in $G(1)$.

At this stage one is in a position to specify the form of the $n$-fold Wild product of $\phi_0$ with itself, corresponding to $\tilde{\gamma} \in G(n)$, indicated with $c_{\tilde{\gamma}}$ in (4):
$$c_{\tilde{\gamma}}(\xi; \phi_0) = \int_{[0, 2\pi)^\infty} \left( \prod_{j \in \tilde{\gamma}} \phi_0(\pi_j \xi) \right) u^\infty(d\theta) \qquad [\tilde{\gamma} \in G(\tilde{\nu}_t), \xi \in \mathbb{R}]. \tag{9}$$

See McKean (1966) and McKean (1967). Then, conditionally on $\tilde{\gamma}$ in $G(\tilde{\nu}_t)$, $c_{\tilde{\gamma}}$ is a mixture, directed by $u^\infty$, of characteristic functions of linear combinations, with coefficients
$$(\pi_1, \ldots, \pi_{\tilde{\nu}_t})(\tilde{\gamma}, \theta),$$

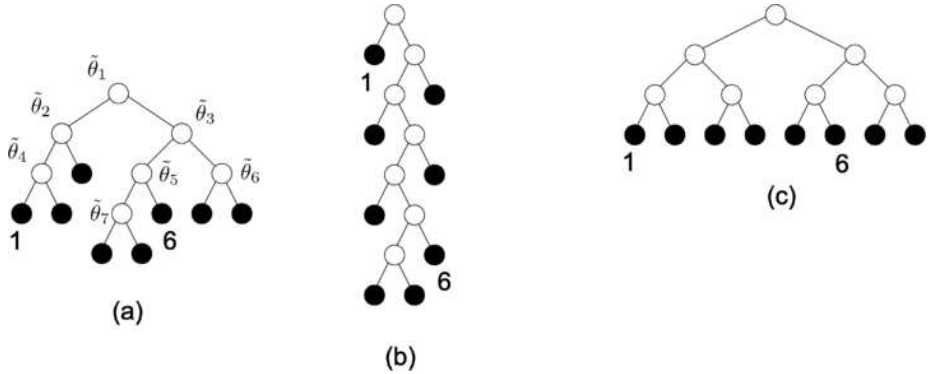

FIG. 1. *Shaded (unshaded) circles stand for leaves (nodes).*

of independent random variables $\tilde{x}_1, \ldots, \tilde{x}_{\tilde{\nu}_t}$. Hence, in view of (3) and (4), one recovers the interpretation of $\mu(\cdot, t)$ stated in Theorem 3 and this completes the proof of the same theorem.

Now we are ready to yield a physical interpretation of this result. If one thinks of each leaf of a tree $\tilde{\gamma}$ with $\tilde{\nu}_t$ leaves as a particle which collides with the particle under observation, then the velocity $V_t$ of this last particle turns out to be the outcome of $\tilde{\nu}_t$ contributions. The $j$th contribution to $V_t$ is given by the initial velocity $\tilde{x}_j$ multiplied by a reducing factor $\pi_j$, which depends on the number $\tilde{\delta}_j$ of collisions that particle $j$ experiences before it collides with the molecule under observation, and on the scattering angles $\theta$. The collisions experienced by particle $j$ take place according to the "order" schematized by $\tilde{\gamma}$.

There is a preliminary statement that plays an important role throughout the rest of the paper. It is drawn from Gabetta and Regazzini (2006) and gives the exact expression of the conditional expectation of $\sum_{j=1}^{\tilde{\nu}_t} x^{\tilde{\delta}_j}$, given $\tilde{\nu}_t$, that is,

$$
(10) \qquad E_t\left(\sum_{j=1}^{\tilde{\nu}_t} x^{\tilde{\delta}_j} \mid \tilde{\nu}_t\right) = \frac{\Gamma(2x + \tilde{\nu}_t - 1)}{\Gamma(2x)\Gamma(\tilde{\nu}_t)}, \qquad x > 0,
$$

which yields

$$
(11) \qquad E_t\left(\sum_{j=1}^{\tilde{\nu}_t} x^{\tilde{\delta}_j}\right) = \sum_{n \geq 1} q_t(n) \frac{\Gamma(2x + n - 1)}{\Gamma(2x)\Gamma(n)} = \exp\{-t(1 - 2x)\}.
$$

Equalities (10)–(11) can be utilized to discuss the asymptotic behavior (as $t \to +\infty$) of the distribution of the random variable

$$
\pi_t^\circ := \max_{1 \leq j \leq \tilde{\nu}_t} |\pi_j|
$$

involved, for example, with the proof of Theorem 2. The starting point for this discussion is given by the following:

LEMMA 1. *For every $x$ in $(0,1)$ and $p > 2$, one has*

$$
P_t\{\pi_t^\circ > x\} \leq \frac{1}{x^p} \exp\{-t(1 - 2\alpha_p)\}.
$$

In particular, $P_t\{\pi_t^\circ > x\} \to 0$, as $t \to +\infty$, even if

$$
(12) \qquad x = x_t = e^{-tc},
$$

provided that $0 < c < (1 - 2\alpha_p)/p$.

For the proof of Lemma 1, see the Appendix.

**3. Proofs of Theorems 1 and 2.** It is useful to premise a remark about the real ($Re$) and imaginary ($Im$) parts of the solution of (2). In fact, it is easy to prove that $Re\phi(\cdot,t)$ is the unique solution of the same problem as (2) with initial data $Re\phi_0$, while $Im\phi(\cdot,t)$ has an explicit form, that is,

$$Im\phi(\xi,t) = (Im\phi_0(\xi))e^{-t}.$$

Then, one can prove Theorems 1 and 2 by assuming, temporarily, that $\phi_0$ is a real-valued characteristic function, which is tantamount to admitting that $\mu_0$ is symmetric, that is,

$$\mu_0((-\infty,-x]) = \mu_0([x,+\infty)) \qquad \text{for every } x > 0.$$

3.1. *Proof of Theorem 2 and of sufficiency in Theorem 1.* We begin with Theorem 2 which, among other things, entails the sufficiency part of Theorem 1. The starting point is an estimate of $|\tilde{\Delta}(\xi)|$, where

$$\tilde{\Delta}(\xi) := \tilde{\phi}_t(\xi) - e^{-\xi^2/2}$$

and $\tilde{\phi}_t$ denotes the conditional characteristic function of

$$\frac{1}{\sigma}\sum_{j=1}^{\tilde{\nu}_t} \pi_j \cdot \tilde{x}_j,$$

given $(\tilde{\nu}_t, \tilde{\gamma}, \tilde{\theta})$:

For every $\varepsilon > 0$ and $\xi$ in $\mathbb{R}$,

$$|\tilde{\Delta}(\xi)| \le e^{-\xi^2/2} \sum_{j=1}^{\tilde{\nu}_t} E_0\bigg[\xi^2 \frac{\pi_j^2 \tilde{x}_j^2}{\sigma^2} \mathbb{1}(|\pi_j \tilde{x}_j| > \sigma\varepsilon)$$

$$+ \varepsilon|\xi|^3 \frac{\pi_j^2 \tilde{x}_j^2}{\sigma^2} \mathbb{1}[|\pi_j \tilde{x}_j| \le \sigma\varepsilon] + \xi^4 \pi_j^2 (\pi_t^\circ)^2\bigg],$$

where $E_0$ indicates expectation with respect to $\mu_0^\infty$ and

$$\varepsilon := (\pi_t^\circ)^a$$

for some $a$ in $(0,1)$, $\pi_t^\circ$ being the same as in Lemma 1. The above inequality follows from a well-known "sharp" estimate of the remainder in the Taylor expansion of $\exp(it)$. A complete proof can be found in Section 9.1 of Chow and Teicher (1997).

Now,

$$\sum_{j=1}^{\tilde{\nu}_t} E_0\bigg[\xi^2 \frac{\pi_j^2 \tilde{x}_j^2}{\sigma^2} \mathbb{1}(|\pi_j \tilde{x}_j| > \sigma\varepsilon)\bigg] \le \sum_{j=1}^{\tilde{\nu}_t} \xi^2 \pi_j^2 E_0\bigg[\mathbb{1}(|\pi_t^\circ \tilde{x}_j| > \sigma\varepsilon)\frac{\tilde{x}_j^2}{\sigma^2}\bigg]$$

$$= \left(\frac{\xi}{\sigma}\right)^2 \int_{|\pi_t^\circ x| > \sigma\varepsilon} x^2 \mu_0(dx)$$

$$\leq \left(\frac{\xi}{\sigma}\right)^2 \int_{|x|>\sigma(\pi_t^\circ)^{a-1}} x^2 \mu_0(dx);$$

$$E_0\left[\sum_{j=1}^{\tilde{\nu}_t} \varepsilon|\xi|^3 \frac{\pi_j^2 \tilde{x}_j^2}{\sigma^2} \mathbb{1}(|\pi_j \tilde{x}_j| \leq \sigma\varepsilon)\right] \leq \sum_{j=1}^{\tilde{\nu}_t} \varepsilon|\xi|^3 \pi_j^2 = |\xi|^3 (\pi_t^\circ)^a$$

and

$$E_0\left[\xi^4 \sum_{j=1}^{\tilde{\nu}_t} \pi_j^2 (\pi_t^\circ)^2\right] = \xi^4 (\pi_t^\circ)^2.$$

Hence,

$$|E_t(e^{i\xi V_t/\sigma}) - e^{-\xi^2/2}| \leq E_t|\tilde{\Delta}(\xi)|$$

(13)
$$\leq \left(\frac{\xi}{\sigma}\right)^2 E_t\left(\int_{|x|>\sigma(\pi_t^\circ)^{a-1}} x^2 \mu_0(dx)\right)$$
$$+ |\xi|^3 E_t((\pi_t^\circ)^a) + \xi^4 E_t((\pi_t^\circ)^2).$$

Next,

$$E_t\left(\int_{|u|>\sigma(\pi_t^\circ)^{a-1}} u^2 \mu_0(du)\right) = E_t\left(\int_{|u|>\sigma(\pi_t^\circ)^{a-1}} u^2 \mu_0(du) \cdot \mathbb{1}\{\pi_t^\circ \leq x\}\right.$$
$$\left.+ \int_{|u|>\sigma(\pi_t^\circ)^{a-1}} u^2 \mu_0(du) \cdot \mathbb{1}\{\pi_t^\circ > x\}\right)$$

$$(x > 0)$$

$$\leq \int_{|u|>\sigma x^{a-1}} u^2 \mu_0(du) + \sigma^2 P\{\pi_t^\circ > x\},$$

which, for $x = x_t := e^{-cpt}$ like in (12), gives

(14) $\quad E_t\left(\int_{|u|>\sigma(\pi_t^\circ)^{a-1}} u^2 \mu_0(du)\right) \leq \int_{|u|>\sigma(x_t)^{a-1}} u^2 \mu_0(du) + \sigma^2 e^{-t(1-2\alpha_p-cp)}.$

Moreover,

(15)
$$E_t((\pi_t^\circ)^a) = E_t((\pi_t^\circ)^a \mathbb{1}\{\pi_t^\circ \leq x_t\}) + E_t((\pi_t^\circ)^a \mathbb{1}\{\pi_t^\circ > x_t\})$$
$$\leq x_t^a + e^{-t(1-2\alpha_p-cp)}$$

and

(16) $$E_t((\pi_t^\circ)^2) \leq x_t^2 + e^{-t(1-2\alpha_p-cp)}.$$

Then, by (13), (14), (15) and (16),

$$|E_t(e^{i\xi V_t/\sigma}) - e^{-\xi^2/2}|$$

$$\leq \left(\frac{\xi}{\sigma}\right)^2 \int_{|u|>\sigma(x_t)^{a-1}} u^2 \mu_0(du)$$
$$+ (|\xi|^2 + |\xi|^3 + |\xi|^4) a_1 e^{-b_1 t} + c_1(|\xi|^3 + |\xi|^4) e^{-c_2 t}$$
$$\leq \left(\frac{\xi}{\sigma}\right)^2 \int_{|u|>\sigma(x_t)^{a-1}} u^2 \mu_0(du) + (\xi^2 + 2|\xi|^3 + 2|\xi|^4) e^{-Bt}$$

holds true with $B = (acp) \wedge (1 - 2ap - cp)$, for every $\xi$. Hence, by Esseen's inequality [see, e.g., Section 9.1 of Chow and Teicher (1997)],

$$\sup_{x\in\mathbb{R}} \left| P_t\left\{\frac{1}{\sigma}V_t \leq x\right\} - G_1(x) \right|$$
$$\leq \frac{2}{\pi} \int_0^T \frac{1}{\xi} \left\{ \left(\frac{\xi}{\sigma}\right)^2 \int_{|u|>\sigma(x_t)^{a-1}} u^2 \mu_0(du) + (\xi^2 + 2|\xi|^3 + 2\xi^4) e^{-Bt} \right\} d\xi$$
$$+ \frac{24}{\sqrt{2\pi^3}} \frac{1}{T}.$$

Then, putting

$$M(t) = \max\left\{ \int_{|u|>\sigma(x_t)^{a-1}} u^2 \mu_0(du), e^{-Bt} \right\}$$

and

$$T = M(t)^{-\beta},$$

one gets

$$\sup_{x\in\mathbb{R}} \left| P_t\left\{\frac{1}{\sigma}V_t \leq x\right\} - G_1(x) \right| \leq B_1 M(t)^{1-4\beta} + B_2 M(t)^{\beta}$$

$$[= BM(t)^{1/5} \quad \text{when } \beta = 1/5].$$

To complete the proof of the first part of Theorem 2, recall that we have stated the previous inequality with initial distribution characterized by $Re(\phi_0)$. Then, for arbitrary characteristic functions $\phi_0$ as initial data, one gets

$$\sup_{x\in\mathbb{R}} |\mu((-\infty,x],t) - G_\sigma(x)| \leq BM(t)^{1/5} + e^{-t}\frac{1}{2}\sup_{x\in\mathbb{R}}|F_0(x) - F_{0,d}(x)|.$$

If $\bar{m}_{2+\delta}$ is finite and $\mu_0$ is symmetric, then from the Berry–Esseen inequality [see, e.g., Theorem 3 in Section 9.1 of Chow and Teicher (1997)],

$$\sup_{x\in\mathbb{R}} |F(x,t) - G_1(x)| \leq \frac{C_\delta}{\sigma^{2+\delta}} E_t\left( \sum_{j=1}^{\tilde{\nu}_t} |\pi_j|^{2+\delta} \overline{m}_{2+\delta} \right)$$
$$\leq \frac{C_\delta}{\sigma^{2+\delta}} \overline{m}_{2+\delta} E_t\left( \sum_{j=1}^{\tilde{\nu}_t} \alpha_{2+\delta}^{\tilde{\delta}_j} \right)$$

$$= \frac{C_\delta}{\sigma^{2+\delta}} \overline{m}_{2+\delta} \sum_{n\geq 1} q_t(n) \frac{\Gamma(2\alpha_{2+\delta}+n-1)}{\Gamma(2\alpha_{2+\delta})\Gamma(n)} \qquad \text{[from (10)]}$$

$$= C_\delta \frac{\overline{m}_{2+\delta}}{\sigma^{2+\delta}} \exp\{-t(1-2\alpha_{2+\delta})\} \qquad \text{[from (11)]}.$$

3.2. *Proof of necessity in Theorem 1.* It remains to prove the *only if* part of Theorem 1. Whence, we assume that the distribution of $V_t$ converges weakly to some probability law on $(\mathbb{R}, \mathscr{B}(\mathbb{R}))$. Moreover, since $\operatorname{Im} \phi(\cdot, t) \to 0$, as $t \to +\infty$, we can confine ourselves to dealing with symmetric initial data, that is, with real-valued $\phi_0$. According to the guidelines indicated at the end of Section 1.2, it is worth recalling that Theorem 3 yields the following representation for the distribution of $V_t$:

$$(17) \qquad P_t\{V_t \in A\} = \int_\Omega \Lambda_{\tilde{\nu}_t}(A, \omega) P_t(d\omega) \qquad [A \in \mathscr{B}(\mathbb{R})],$$

where $\Lambda_{\tilde{\nu}_t}$ indicates the $\tilde{\nu}_t$-fold convolution of $\lambda_{1,t}, \ldots, \lambda_{\tilde{\nu}_t,t}$, $\lambda_{j,t}$ standing for a *conditional distribution* of $\pi_j \tilde{x}_j$, given $(\tilde{\nu}_t, \tilde{\gamma}, \tilde{\theta})$, for $j = 1, \ldots, \tilde{\nu}_t$. Now, following the above guidelines, let us analyze the asymptotic behavior (as $t \to +\infty$) of $\Lambda_{\tilde{\nu}_t,t}$ together with that of all the elements which figure in general formulations of the central limit theorem, that is,

$$W_t := (\Lambda_{\tilde{\nu}_t,t}, \lambda_{1,t}, \ldots, \lambda_{\tilde{\nu}_t,t}, \delta_0, \ldots, \tilde{\gamma}, \tilde{\theta}, \tilde{\nu}_t, U_t(1/2), U_t(1/3), \ldots),$$

where $\delta_y$ stands for the unit mass at $y$ and, for any $\zeta > 0$,

$$U_t(\zeta) := \operatorname{Max}_{1\leq j\leq \tilde{\nu}_t} \lambda_{j,t}([-\zeta, \zeta]^c).$$

To grasp the importance of the elements of $W_t$, it is worth recalling the classical formulation of the central limit theorem for independent *uniformly asymptotically negligible* (*uan*) summands $X_{nk}$ ($k = 1, \ldots, m_n$, $n = 1, 2, \ldots$) with symmetric distributions ($F_{nk}$ will denote the probability distribution function of $X_{nk}$ for every $k$ and $n$):

*In order that $\sum_{k=1}^{m_n} X_{nk}$ can converge in distribution, it is necessary and sufficient that there exist a nonnegative number $\sigma$ and a symmetric Lévy measure $l$ (a measure on $\mathbb{R}\setminus\{0\}$ satisfying $\int_{\mathbb{R}\setminus\{0\}} (y^2 \wedge 1) l(dy) < +\infty$ and $l((-\infty, -x]) = l([x, +\infty))$ for every $x > 0$) which meets the following conditions:*

$$(18) \quad l([x, +\infty)) = \lim_{n\to +\infty} \sum_k \{1 - F_{nk}(x)\} \qquad \text{if } x > 0 \text{ and } l\{x\} = 0$$

$$(19) \qquad \sigma^2 = \lim_{\varepsilon\to 0^+} \overline{\lim}_n \sum_k \int_{[-\varepsilon, \varepsilon]} x^2 \, dF_{nk}(x).$$

*In case these conditions are satisfied, the limiting distribution of $\sum_k X_{nk}$ is the infinitely divisible law characterized by the Fourier–Stieltjes transform*

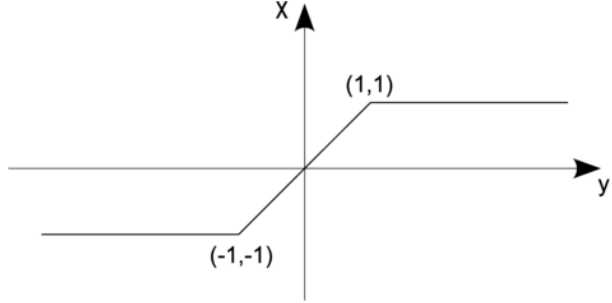

Fig. 2.  *The function $\chi$.*

$\exp\{-\psi\}$ *with*

(20) $$\psi(u) = \frac{\sigma^2 u^2}{2} + \int_{\mathbb{R}\setminus\{0\}} (1 - e^{-iuy} + iu\chi(y))l(dy),$$

$\chi$ *being the function shown in Figure* 2.

This specific version of the central limit theorem is drawn from Section 16.8 of Fristedt and Gray (1997).

Think of the range of $W_t$ as a subset of

$$S := \mathbb{P}(\overline{\mathbb{R}})^\infty \times G^* \times [0, 2\pi)^\infty \times \overline{\mathbb{R}}^\infty,$$

where, given any metric space $M$, $\mathbb{P}(M)$ stands for the set of all probability measures on the Borel class $\mathscr{B}(M)$; $\overline{\mathbb{R}} = [-\infty, +\infty]$ is equipped with the distance $d(x,y) := |\arctan y - \arctan x|$ for any $(x,y) \in \overline{\mathbb{R}}^2$. It is well known that $\mathbb{P}(\overline{\mathbb{R}})$ can be metrized, consinstenly with the topology of weak convergence, in such a way that it may turn out to be a compact and separable metric space; see Sections 5 and 6 (vi) of Billingsley (1999). Moreover, think of the set $G$ of all McKean's trees as a metric space with the discrete distance, and define $G^*$ to be a metrizable compactification of $G$, which exists since $G$ is separable; see, for example, Corollary 1 in Section 10.1 of Gemignani (1990). Therefore, $S$ proves to be a separable and compact metric space with respect to the product topology. Hence, any family of probability measures on $(S, \mathscr{B}(S))$ is tight; in particular, the family of the probability distributions $Q_t$ of $W_t$, $t > 0$, turns out to be tight. At this stage, the conclusive steps of the proof rest on the following lemmata.

LEMMA 2. *For every positive $\delta$ and $\beta$, one has*

$$P\{U_t(\delta) > \beta\} \to 0 \qquad (t \to +\infty).$$

LEMMA 3. *If the law of $V_t$ converges weakly as $t \to +\infty$, then any sequence $(Q_{t_n})_n$ of elements of $\{Q_t : t > 0\}$, such that $t_n \nearrow +\infty$, contains a* subsequence $(Q_{t_{n'}})_{n'}$ weakly convergent to a probability measure $Q$ supported by

$$\mathbb{P}_0(\overline{\mathbb{R}}) \times \{\delta_0\}^\infty \times G^* \times [0, 2\pi)^\infty \times \{+\infty\} \times \{0\}^\infty$$

with

$$\mathbb{P}_0(\overline{\mathbb{R}}) := \{p \in \mathbb{P}(\overline{\mathbb{R}}) : p(\{-\infty, +\infty\}) = 0\}.$$

Whence, since $S$ is separable, from the Dudley generalization of a Skorokhod's theorem [see, e.g., Theorem 11.7.2 in Dudley (2002)], one can apply Lemmatas 1, 2 and 3 to state that, on some probability space $(\Omega^*, \mathscr{F}^*, P^*)$ there are random elements

$$W^*_{t_{n'}} = (\Lambda^*_{\tilde{\nu}^*_{t_{n'}}}, \lambda^*_{1,t_{n'}}, \ldots, \lambda^*_{\tilde{\nu},t_{n'}}, \delta_0, \ldots, \tilde{\gamma}^*, \tilde{\theta}^*, \tilde{\nu}^*_{t_{n'}}, U^*_{t_{n'}}(1/2), \ldots)$$

taking values in $S$, with distribution $Q_{t_{n'}}$, satisfying

$$\Lambda^*_{\tilde{\nu}^*_{t_{n'}}} \Rightarrow \Lambda^*, \qquad \lambda^*_{j,t_{n'}} \Rightarrow \delta_0 \qquad (j = 1, 2, \ldots),$$

(21)

$$\tilde{\nu}^*_{t_{n'}} \to +\infty, \qquad U^*_{t_{n'}}(1/k) \to 0$$

for $k = 2, 3, \ldots$ on a set $\Omega_1^*$ of $\mathscr{F}^*$ such that $P^*(\Omega_1^*) = 1$, provided that $(t_{n'})$ is the same subsequence $(t_{n'})$ as in Lemma 3. (The symbol $\Rightarrow$ is used to designate weak convergence of probability measures.) The distributional properties of $W^*_{t_{n'}}$ imply that $\Lambda^*_{\tilde{\nu}^*_{t_{n'}}}$ is the convolution of $\lambda^*_{1,t_{n'}}, \ldots, \lambda^*_{\tilde{\nu}^*_{t_{n'}},t_{n'}}$, and that equality $U^*_{t_{n'}}(1/k) = \text{Max}_{1 \leq j \leq \tilde{\nu}^*_{t_{n'}}} \lambda^*_{j,t_{n'}}([-1/k, 1/k]^c)$ holds true for every $k$. Thus, conditions (18)–(19) must be valid with $\lambda^*_{j,t_{n'}}((-\infty, \cdot])$ in the place of $F_{t_{n'},j}(\cdot)$. Apropos of (19), note that

(22)
$$\sum_{j=1}^{\tilde{\nu}^*_{t_{n'}}} \int_{|x|<\varepsilon} x^2 \lambda^*_{j,t_{n'}}(dx) \geq \sum_{j=1}^{\tilde{\nu}^*_{t_{n'}}} (\pi^*_j)^2 \int_{\{x:\,|\pi^*_j x|<\varepsilon\}} x^2 \mu_0(dx)$$
$$\geq \int_{\{x:\,|x|(\pi^\circ_{t_{n'}})^*<\varepsilon\}} x^2 \mu_0(dx),$$

with $(\pi^\circ_t)^* = \max_{1 \leq j \leq \tilde{\nu}^*_t} |\pi^*_j|$. From Lemma 1, combined with a well-known necessary and sufficient condition, for convergence in probability, in terms of sub–subsequences converging almost surely [see, e.g., Lemma 2 in Section 3.3 of Chow and Teicher (1997)], there is a subsequence $(t_{n''})$ of $(t_{n'})$ such that $(\pi^\circ_{t''})^* \to 0$ ($P^*$-almost surely). Hence, from (19) and (22), it turns out that $\int_\mathbb{R} x^2 \mu_0(dx) = \lim_{t_{n''}} \int_{\{x:|x|(\pi^\circ_{t_{n''}})^*<\varepsilon\}} x^2 \mu_0(dx)$ must be finite. This completes the proof of Theorem 1 when $\mu_0$ is symmetric. The extension to general initial data follows from the simple remark that the second moment of $\mu_0$ is finite if and only if the second moment of the "even" component of $\mu_0$ is finite.

## APPENDIX

PROOF OF LEMMA 1. For any $x > 0$,

$$\begin{aligned} P_t\{\pi_t^\circ \leq x\} &= 1 - P_t\left(\bigcup_j \{|\pi_j| > x\}\right) \\ &\geq 1 - \sum_j P_t\{|\pi_j| > x\} \\ &\geq 1 - \frac{1}{x^p} \sum_j E_t(|\pi_j|^p) \quad \text{(from the Markov inequality)} \\ &= 1 - \frac{1}{x^p} \sum_j E_t\left(\prod_{i=1}^{\delta_j} |\tau_i^{(j)}|^p\right) \\ &= 1 - \frac{1}{x^p} \sum_j E_t(\alpha_p^{\delta_j}) \quad \text{(from Theorem 3)} \\ &= 1 - \frac{1}{x^p} \exp\{-t(1 - 2\alpha_p)\} \quad \text{[from (11)]}. \quad \square \end{aligned}$$

PROOF OF LEMMA 2. Fix $\beta > 0$ and sufficiently small $\varepsilon$ so that

$$\mu_0(\{x : |x| > \delta/\varepsilon\}) \leq \beta.$$

Now, notice that

$$\begin{aligned} P_t\{U_t(\delta) > \beta\} &\leq P_t\{\pi_t^\circ > \varepsilon\} \\ &\quad + P_t\{\pi_t^\circ \leq \varepsilon, \mu_0(\{x : |x| > \delta/\varepsilon\}) > \beta\} \\ &= P_t\{\pi_t^\circ > \varepsilon\} \end{aligned}$$

and apply Lemma 1. $\square$

PROOF OF LEMMA 3 [FROM FORTINI, LADELLI AND REGAZZINI (1996)]. In view of the tightness of $\{Q_t : t > 0\}$, the Prokhorov theorem [cf., e.g., Section 5 of Billingsley (1999)] can be applied to state the existence of a weakly convergent subsequence $(Q_{t_{n'}})$ of $(Q_{t_n})$. From Lemma 2 and the fact that $P_t\{\tilde{\nu}_t > K\} \to 1$, as $t \to +\infty$, for every $K > 0$, it is easy to check that the limiting distribution of

$$(\lambda_{1,t_{n'}}, \ldots, \lambda_{\tilde{\nu}_{t_{n'}}}, \delta_0, \ldots, \tilde{\gamma}, \tilde{\theta}, \tilde{\nu}_{t_{n'}}, U_{t_{n'}}(1/2), \ldots)$$

is supported by $\{\delta_0\}^\infty \times G^* \times [0, 2\pi]^\infty \times \{+\infty\} \times \{0\}^\infty$. Then, it is enough to prove that the weak limit $Q^{(1)}$ of the law $Q_{t_{n'}}^{(1)}$ of $\Lambda_{\tilde{\nu}_{t_{n'}}^*}$ is supported by $\mathbb{P}_0(\overline{\mathbb{R}})$. Since $(V_{t_n})$ converges in law, from a theorem of Le Cam, it must be tight; see Section 5 of Billingsley (1999). Then, for every positive integer $m$, there is $K_m$ satisfying $K_m \nearrow +\infty$ and

$$P_t\{|V_{t_n}| > K_m\} \leq 1/m \qquad (m = 1, 2, \ldots).$$

Now fix $\eta$ in $(0,1)$ and put $[-K_m, K_m]^c = \overline{\mathbb{R}} \setminus [-K_m, K_m]$. Then, from (17),

$$P_t\{|V_{t_n}| > K_m\} = \int \Lambda_{\tilde{\nu}_{t_n}}([-K_m, K_m]^c) \, dP_t$$

$$\geq \int \Lambda_{\tilde{\nu}_{t_n}}([-K_m, K_m]^c) \mathbb{1}_{(\eta, +\infty)}(\Lambda_{\tilde{\nu}_{t_n}}([-K_m, K_m]^c)) \, dP_t$$

$$\geq \eta Q_{t_n}^{(1)}(A_\eta^{(m)}),$$

with

$$A_\eta^{(m)} = \{p : p \in \mathbb{P}(\overline{\mathbb{R}}), p([-K_m, K_m]^c) > \eta\}.$$

Then

$$Q_{t_{n'}}^{(1)}(A_\eta^{(m)}) \leq \frac{1}{m\eta}.$$

A direct application of the Alexandroff "portmanteau" theorem [see, e.g., Theorem 2.1 in Billingsley (1999)] shows that $C_\eta^{(m)} := (A_\eta^{(m)})^c$ is closed. Then, from the same theorem [see point (iii) in Billingsley (1999)] one deduces

$$Q^{(1)}(C_\eta^{(m)}) \geq \overline{\lim_{n'}} Q_{t_{n'}}^{(1)}(C_\eta^{(m)}) \geq 1 - \frac{1}{m\eta}.$$

Clearly, as $m \to +\infty$,

$$Q^{(1)}(C_\eta^{(m)}) \uparrow Q^{(1)}\left(\bigcup_m C_\eta^{(m)}\right) \subset Q^{(1)}(C_\eta^{(\infty)}),$$

with

$$C_\eta^{(\infty)} = \{p : p\{-\infty, +\infty\} \leq \eta\}.$$

Whence, $Q^{(1)}(C_\eta^{(\infty)}) \geq 1 - \frac{1}{m\eta}$ for every $m$ and this entails $Q^{(1)}(C_\eta^{(\infty)}) = 1$ for every $\eta > 0$, which is tantamount to saying that $\{p : p\{-\infty, +\infty\} = 0\}$ has probability one. $\square$

**Acknowledgments.** The authors would like to thank E. Carlen, M. C. Carvalho for advice and encouragement and F. Bassetti for his significant contributions to improve parts of the first draft. They also thank the reviewers for their time, comments and suggestions.

Dipartimento di Matematica
Università degli Studi di Pavia
via Ferrata 1
27100 Pavia
and
IMATI-CNR
Italy
E-mail: ester.gabetta@unipv.it
       eugenio.regazzini@unipv.it